\documentclass[12pt]{amsart}
\usepackage{amssymb}
\usepackage{amsfonts}
\usepackage{subfigure}
\usepackage{graphicx,color}

\newcommand{\abs}[1]{\ensuremath{\left| #1 \right| }}

\newcommand{\ST}{\widehat{S}_{(K)}}
\newcommand{\trace}{\mathrm{trace}}

\theoremstyle{plain}

\newtheorem{lemma}{Lemma}

\newtheorem{remark}{Remark}

\newtheorem{theorem}{Theorem}
\numberwithin{equation}{section}

\begin{document}
\title{The bias-variance trade-off in Thomson's multitaper estimator}
\author{Lu\'{\i}s Daniel Abreu}
\address{Acoustics Research Institute, Austrian Academy of Science,
Wohllebengasse 12-14 A-1040, Vienna Austria}
\email{labreu@kfs.oeaw.ac.at}
\author{Jos\'{e} Luis Romero}
\address{Faculty of Mathematics \\
University of Vienna \\
Oskar-Morgenstern-Platz 1 \\
A-1090 Vienna, Austria}
\email{jose.luis.romero@univie.ac.at}
\subjclass{}
\keywords{}
\thanks{L. D. A. was supported by the Austrian Science Fund (FWF)
START-project FLAME ("Frames and Linear Operators for Acoustical Modeling
and Parameter Estimation") 551-N13. J. L. R. gratefully acknowledges support
from a Marie Curie fellowship, within the 7th. European
Community Framework program, under grant PIIF-GA-2012-327063.}

\begin{abstract}
At the heart of non-parametric spectral estimation, lies the dilemma known
as the \emph{bias-variance trade-off}: low-biased estimators tend to have
high variance and low variance estimators tend to have high bias. In 1982,
Thomson introduced a multitaper method where this trade-off is made explicit
by choosing a target bias resolution and obtaining a corresponding variance reduction.
The method became the standard in many applications. Its favorable
bias-variance trade-off is due to an empirical fact, conjectured by
Thomson based on numerical evidence: assuming bandwidth $W$ and $N$ time domain observations,
\emph{\ the average of the square of the first }$K=\left\lfloor
2NW\right\rfloor $\emph{\ Slepian functions approaches, as }$K$\emph{\
grows, an ideal band-pass kernel for the interval }$\left[ -W,W\right] $. We
provide an analytic proof of this fact and quantify the approximation error
in the $L^{1}$ norm; the approximation error is then used to control the
bias of the multitaper estimator resulting from spectral leakage. This
leads to new performance bounds for the method, explicit \ in terms of the
bandwidth $W$ and the number $N$ of time domain observations. Our method is
flexible and can be extended to higher dimensions and different geometries.
\end{abstract}

\maketitle

\section{Introduction}

Let $I=\left[ -1/2,1/2\right] $. Any stationary, real, ergodic,
zero-mean, Gaussian stochastic process has a \emph{Cram\'{e}r spectral
representation} 
\begin{equation*}
x(t)=\int_{I}e^{2\pi i\xi t}dZ(\xi )\text{,}
\end{equation*}
and the \emph{spectrum} $S(\xi )$, defined as 
\begin{equation*}
S(\xi )d\xi =\mathbb{E}\{\left\vert dZ(\xi )\right\vert ^{2}\}\text{,}
\end{equation*}
and often called the \emph{power spectral density} of the process, yields
the periodic components of $x(t)$. The goal of spectral estimation is to
solve the highly underdetermined problem of \emph{estimating }$S(\xi )$
\emph{\ from a sample of }$N$\emph{\ contiguous observations }$x(0),...,x(N-1)$.
Embryonic approaches to the problem (Stokes 1879, Shuster 1898) used the so
called \emph{periodogram}:
\begin{equation}
\widehat{S}(\xi )=\frac{1}{N}\left\vert \sum_{t=0}^{N-1}x(t)e^{-2\pi i\xi
t}\right\vert ^{2},
\end{equation}
whose analysis has influenced harmonic analysts since
Norbert Wiener (see \cite{Benedetto}). The periodogram can also be weighted with a data window $\left\{ D_{t}\right\}
_{t=0}^{N-1}$, usually called a\emph{\ taper}, giving the estimator:
\begin{equation}
\widehat{S}_{D}(\xi )=\left\vert \sum_{t=0}^{N-1}x(t)D_{t}e^{-2\pi i\xi
t}\right\vert ^{2}\text{.}  \label{direct}
\end{equation}
The choice of the taper $\left\{ D_{t}\right\} _{t=0}^{N-1}$ can have a significant effect on the
resulting spectrum estimate $\widehat{S}_{D}$. This is apparent by observing that its
expectation is the convolution of the \emph{true (nonobservable) spectrum}
$S(\xi )$ with the \emph{spectral window} $\abs{\mathcal{F}D(\xi )}^2=
\abs{\sum_{t=0}^{N-1}D_{t}e^{-2\pi i\xi t}}^2$, i.e., 
\begin{equation}
\label{smooth}
\mathbb{E}\left\{ \widehat{S}_{D}(\xi )\right\} =S(\xi )\ast 
\abs{\mathcal{F}D(\xi )}^2.
\end{equation}
Thus, the bias of the tapered estimator, which is the difference
$S(\xi) - \mathbb{E}\{\widehat{S}_{D}(\xi )\}$, is determined by the
smoothing effect of $\left\{ D_t \right\}_{t=0}^{N-1}$ over the true
spectrum. Ideally, the function
$\mathcal{F}D(\xi )$ should be concentrated on the interval 
$[-\tfrac{1}{2N},\tfrac{1}{2N}]$, but the 
\emph{uncertainty principle of Fourier analysis} precludes such perfect concentration.
Inevitably, some portion of the filter $\mathcal{F}D(\xi )$ will lie outside the target
region and \emph{spectral leakage} occurs.

In \cite{Thomson}, Thomson used the sequences which minimize 
spectral leakage to construct an algorithm using several tapered
estimates, whence the name \emph{multitaper}. In doing so, he was able to
reduce variance by averaging, while introducing a tolerable amount of
spectral leakage. Thomson's multitaper method has been used in a variety of scientific
applications including climate analysis (see, for instance \cite{Science1997},
or \cite{PNAS2012} for a local spherical approach), and it was used to better understand the relation between 
atmospheric $CO_{2}$ and climate change (see \cite[Section 1]{Thomson2}). The method became also paramount
in statistical signal analysis \cite{PW}.

Today, Thomson's multitaper method remains an effective spectral estimation method.
It has recently found remarkable applications in electroencephalography \cite
{Neurosciences} and it is the preferred spectral sensing procedure \cite
{Cognitive} for the rapidly emerging field of cognitive radio \cite{Haykin}.
In the next paragraph we provide an outline of the essence of the method.

Thomson's method starts by selecting a target frequency smoothing band $[-W,W]$ with
$1/2N<W<1/2$, thus accepting a reduction in spectral resolution by a factor of about $2NW$. The \emph{first step} consists of obtaining a number $K=\left\lfloor 2NW\right\rfloor $ (the smallest integer not greater than
$2NW $) of estimates of the form \eqref{direct} by setting, for every $k\in
\{0,...,K-1\}$, $D_{t}=v_{t}^{(k)}(N,W)$, where the \emph{discrete prolate
spheroidal sequences} $v_{t}^{(k)}(N,W)$ are defined as the solutions of the
Toeplitz matrix eigenvalue equation
\begin{equation*}
\sum_{n=0}^{N-1}\frac{\sin 2\pi W\left( t-n\right) }{\pi \left( t-n\right) }
v_{n}^{(k)}(N,W)=\lambda _{k}(N,W)v_{t}^{(k)}(N,W)\text{.}
\end{equation*}
The resulting tapered periodogram is then denoted by $\widehat{S}_{k}(\xi )$. The \emph{second step} consists of averaging. One uses the estimator
\begin{equation}
\ST(\xi )=\frac{1}{K}\sum_{k=0}^{K-1}\widehat{S}_{k}(\xi ),
\label{Thomson}
\end{equation}
which achieves a reduced variance (see \cite{Thomson} for an asymptotic analysis of slowly varying spectra and 
\cite{WMP,liro08} for non-asymptotic expressions).

\begin{figure}
\centering
\subfigure[Slepians $U_k$ and their squares $\abs{U_k}^2$, for $k=1,5,9,19$.]{
\includegraphics[scale=0.23]{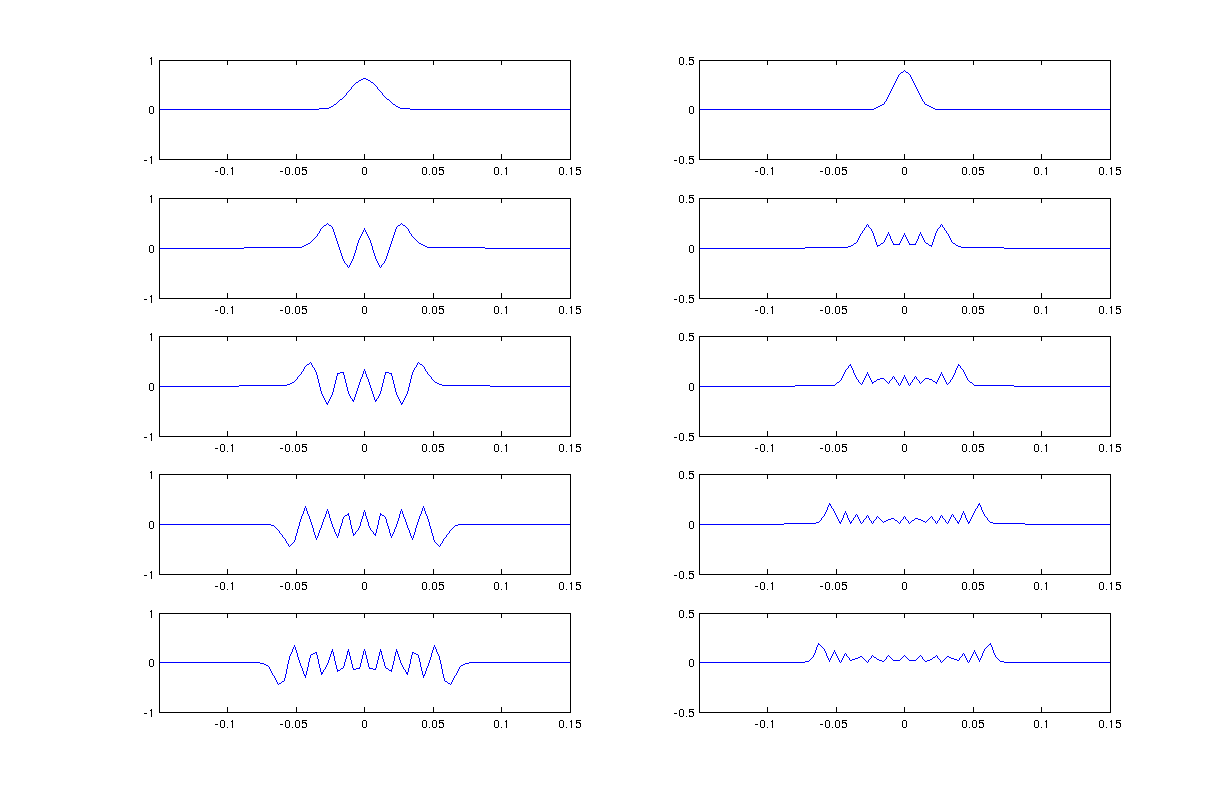}
\label{fig_prolates}
}
\subfigure[Thomson's spectral window.]
{
\includegraphics[scale=0.20]{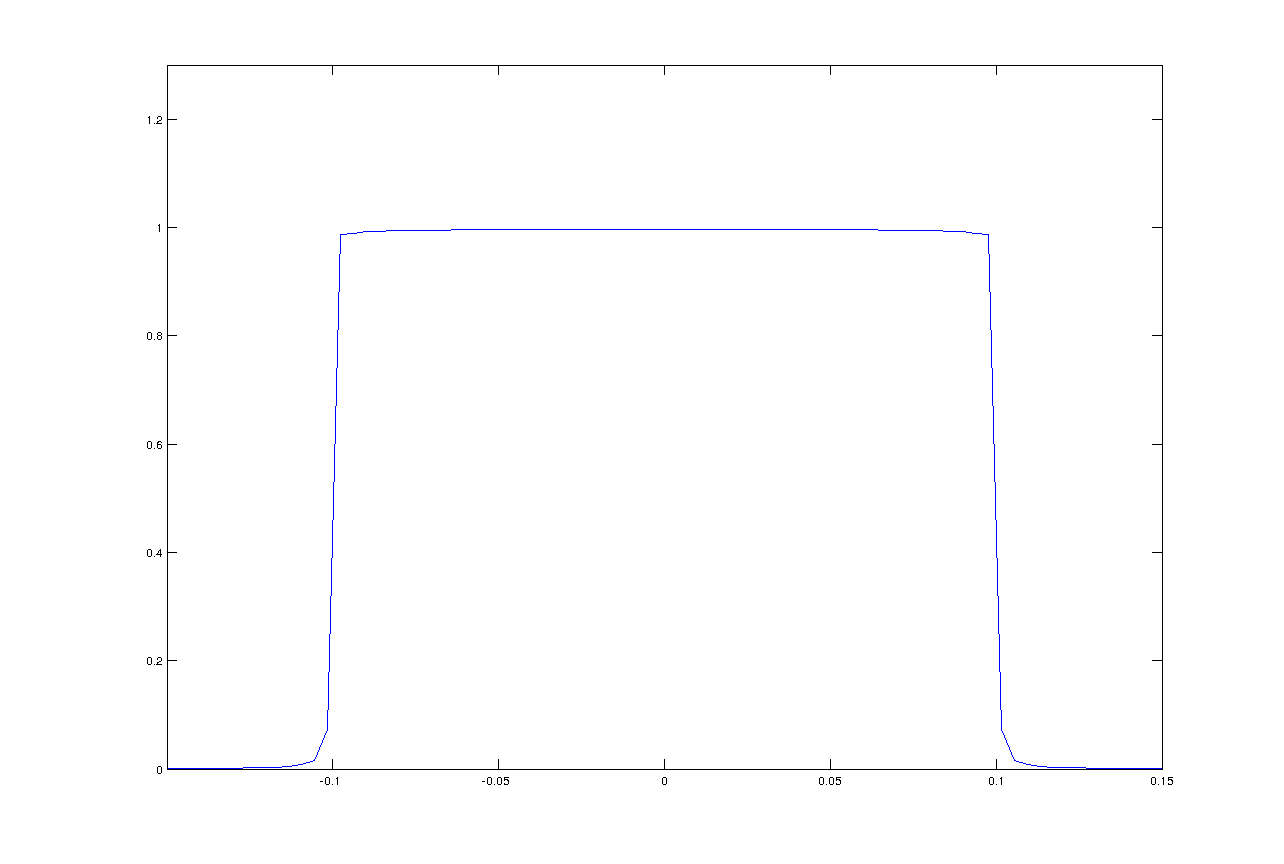}
}
\caption{Some Slepians and the spectral window with $N=256$ and $w=0.1$.}
\label{fig}
\end{figure}

To inspect the performance of the estimator $\ST(\xi )$ on the spectral domain, let us consider the \emph{discrete prolate spheroidal functions}, also known as \emph{Slepians}. They
are the discrete Fourier transforms of the sequences $v_{t}^{(k)}(N,W)$,
denoted by $U_{k}(N,W;\xi )$, and satisfy the integral equation 
\begin{equation}
\int_{-W}^{W}\mathbf{D}_{N}(\xi -\xi {\acute{}})U_{k}(N,W;\xi {\acute{}}
)d\xi {\acute{}}=\lambda _{k}(N,W)U_{k}(N,W;\xi )\text{,}
\label{eq_prolates}
\end{equation}
where
\begin{equation}
\mathbf{D}_{N}(x)=\frac{\sin N\pi x}{\sin \pi x}  \label{Dirichlet}
\end{equation}
is the Dirichlet kernel. Observe that, according to (\ref{smooth}),
$\mathbb{E}\{\widehat{S}_{k}(\xi )\}$ is a smoothing average of the unobservable spectrum by the kernel $\left\vert 
U_{k}(N,W;\xi )\right\vert ^{2}$. Recall that
the bias of each individual estimate in (\ref{Thomson}) is given by 
\begin{equation}
Bias\left( \widehat{S}_{k}(\xi )\right) =\mathbb{E}\{\widehat{S}_{k}(\xi )\}-S(\xi )=S(\xi )\ast\left\vert U_{k}(N,W;\xi 
)\right\vert ^{2}-S(\xi ).  \label{bias}
\end{equation}
The optimal concentration of the first prolate function on the interval
$[-W,W]$ leads to a low bias when $k=0$. But since the amount of energy of
$U_{k}(N,W;\xi )$ inside $[-W,W]$ decreases with $k$ (because the energy is
given by the eigenvalues in (\ref{eq_prolates}) and they decrease from $1$
to $0$ as $k$ approaches $K$), the bias increases with $k$. 
To explain the remarkable performance of the \emph{averaged} estimator, Thomson noted the following: the expected value of the estimator (\ref{Thomson}) is given
by 
\begin{equation}
\mathbb{E}\{\ST(\xi )\}=\frac{1}{K}\sum_{k=0}^{K-1}\mathbb{E}\{\widehat{S}_{k}(\xi )\}=S(\xi ) \ast \frac{1}{K}\rho _{K}(N,W;\xi )\text{,}  \label{Sk}
\end{equation}
where 
\begin{equation}
\frac{1}{K}\rho _{K}(N,W;\xi )=\frac{1}{K}\sum_{k=0}^{K-1}\left\vert
U_{k}(N,W;\xi )\right\vert ^{2}  \label{eq_inten}
\end{equation}
is the \emph{spectral window} of (\ref{Sk}). The bias performance is due to the fact,
numerically illustrated by Thomson, that the spectral window (\ref{eq_inten}) is very similar
to a flat function localized on $[-W,W]$ (see Figure \ref{fig}). This is an intriguing mathematical
phenomenon. Heuristically, it requires the functions in the sequence
$\{\left\vert U_{k}(N,W,\cdot )\right\vert ^{2}:k=0,\ldots, K-1\}$ to be
organized inside the interval $[-W,W]$ in a very particular way: \emph{each
function tends to fill in the empty energy spots left by the sum of the
previous ones}. This behavior is reminiscent of the Pythagorean relation for
pure frequencies: $\sin ^{2}(t)+\cos ^{2}(t)=1$. More precisely, claiming
that the spectral window in Thomson's method approximates an ideal band-pass
kernel, means that the two functions 
\begin{equation}
\frac{1}{K}\rho _{K}(N,W,.)\text{ \ \ \ and \ \ \ \ }\frac{1}{2W}\mathbf{1}_{[-W,W]}\text{,}  \label{two}
\end{equation}
approach each other as $K$ increases. This is indeed true and we provide an
analytic bound for the $L^{1}$-distance between the functions in \eqref{two}.

\begin{theorem}[Spectral leakage estimate]
\label{th_main} Let $N\geq 2$ be an integer, $W\in (-1/2,1/2)$ and set $K:=\left\lfloor 2NW\right\rfloor $. Then 
\begin{equation}
\left\Vert \frac{1}{K}\rho _{K}(N,W,\cdot )-\frac{1}{2W}\mathbf{1}_{[-W,W]}\right\Vert _{L^{1}(I)}\lesssim \frac{\log N}{K}.
\label{L1}
\end{equation}
\end{theorem}

The spectral leakage estimate (\ref{L1}) is precisely what we need in order to 
quantify Thomson's asymptotic analysis of the bias
of the multitaper estimator \cite[pag. 1062]{Thomson}
and validate the bias-variance trade-off. This is explained in
the Conclusion section.

A relevant feature of the method introduced in this paper is its
flexibility. While the description of each individual solution
to the concentration problem in \eqref{eq_prolates} is very subtle, the \emph{aggregated behavior} of the critical number of solutions to \eqref{eq_prolates} displays a simple profile. A similar aggregated behavior has been investigated in \cite{AGR} and numerically illustrated in \cite{BB, SimmonsGeoph2008}.

Our analysis depends on the properties of the eigenvalues in \eqref{eq_prolates}.
Similar properties have been recognized in the
eigenvalue problem in the context of Hankel bandlimited
functions \cite{AB}. Since the problem studied in \cite{AB} includes the one considered by Slepian in his construction of $2d$ radial prolate functions, we expect our methodology to be applicable to spectral estimation problems involving $2d$ functions whose spectrum lies on a disk. This may have applications in cryo-electron microscopy, where estimation of noise stochastics is an important
consideration when applying PCA to microscopy images \cite{zhsi13}.
Other multitaper estimators include multi-window estimators for
non-stationary spectrum \cite{BB,OW} and the one based on spherical Slepians \cite{SIAMRev,SimmonsGeoph2008}.

\section{Proof of the main result}
Our proof uses tools from the Landau-Pollack-Slepian theory
\cite{Slepian,posl61,la67-1,la75-1,lapo61}. We do not rely on special properties of the interval $I$, but rather on 
so-called trace / norm estimates that can be obtained in many other contexts of practical interest (e.g. \cite{AB}). 
Hence the flexibility of our approach.

Let $I:=[-1/2,1/2]$ and let us denote the exponentials by $e_{\omega
}(x):=e^{2\pi ix\omega }$. We will always let $N\geq 2$ be an integer and $W\in (-1/2,1/2)$. For two
non-negative functions $f,g$, the notation $f\lesssim g$ means that there
exists a constant $C>0$ such that $f\leq Cg$. (The constant $C$, of course,
does not depend on the parameters $N,W$.)

\subsection{Trigonometric polynomials}

For notational convenience, we use a temporal normalization that is slightly
different of the one in the Introduction (this has no impact in the
announced estimates). We consider the space of trigonometric polynomials 
\begin{equation*}
\mathcal{P}_{N}=Span\left\{ e_{\frac{-N+1}{2}+j}:0\leq j\leq N-1\right\}
\subseteq L^{2}\left( I\right) \text{.}
\end{equation*}
This is a Hilbert space with a reproducing kernel given by the translated
Dirichlet kernel, $\mathbf{D}_{N}(x-y)$, \ $x,y\in I$,\ $N\in \mathbb{N}$,
where $\mathbf{D}_{N}$ is given by (\ref{Dirichlet}). Note that $\int_{I}\left\vert \mathbf{D}_{N}\right\vert ^{2}=N$.

\subsection{Toeplitz operators}
For $W\in (-1/2,1/2)$ the \emph{Toeplitz
operator} $H_{W}^{N}$ is 
\begin{equation}
\label{eq_toep_op}
H_{W}^{N}f:=P_{{\mathcal{P}_{N}}}\left( (P_{\mathcal{P}_{N}}f)\cdot
1_{[-W,W]}\right) ,\qquad f\in L^{2}(I),
\end{equation}
where $P_{{\mathcal{P}_{N}}}$ is the orthogonal projection onto $\mathcal{P}_N$.
When $f\in \mathcal{P}_{N}$, $H_{W}^{N}f$ is simply the projection of $f\cdot 1_{[-W,W]}$ into $\mathcal{P}_{N}$. The 
Slepian functions $\{U_{k}(N,W) : k=0, \ldots, N-1\}$
are the eigenfunctions of $H_{W}^{N}$ with corresponding eigenvalues $\lambda_k=\lambda_K(N,W)$:
\begin{align}
\label{eq_eigen}
\int_{-W}^W \abs{U_{k}(N,W;\xi)}^2 \, d\xi = \lambda_k,
\end{align}
ordered non-increasingly. We normalize the Slepian functions by:
$\int_I \abs{U_{k}(N,W;\xi)}^2 \, d\xi = 1$.
We will need a description of the profile of the eigenvalues of $H^N_W$.
\begin{lemma}
\label{lemma_sum_eig}
For $N \geq 2$, $W \in (-1/2,1/2)$ and $K:=\left\lfloor 2NW\right\rfloor$:
\begin{equation}
\label{eq_bound}
\abs{1 - \frac{1}{K}\sum_{k=0}^{K-1} \lambda_k(N,W)} \lesssim \frac{\log N}{K}.
\end{equation}
\end{lemma}
We postpone the proof of Lemma \ref{lemma_sum_eig} to the Appendix.
The quantity on the left-hand side of \eqref{eq_bound} has been studied in 
\cite{liro08} to qualitatively analyze the performance of Thomson's method. Lemma \ref{lemma_sum_eig} 
refines the analysis of \cite{liro08}, giving a concrete growth estimate. (See also the remarks 
after Theorem 5 in \cite{liro08}.)

\subsection{\textbf{Proof of Theorem }\protect\ref{th_main}}
We first estimate the narrow band error. Note that $\rho _{K}(N,W;\xi) = \sum_{k=0}^{K-1} \abs{U_{k}(N,W;\xi)}^2
\leq \sum_{k=0}^{N-1} \abs{U_{k}(N,W;\xi)}^2 = D_N(0) = N$. Consequently,
$\tfrac{1}{K}\rho_K(N,W;\xi) \leq \tfrac{N}{K}$ and, using \eqref{eq_eigen}, we can estimate:
\begin{align*}
&\int_{-W}^W \abs{\frac{1}{K}\rho _{K}(N,W;\xi)-\frac{1}{2W}\mathbf{1}_{[-W,W]}(\xi)} \,d\xi
\\
&\qquad\leq
\int_{-W}^W \abs{\left(\frac{1}{2W}-\frac{N}{K}\right)\mathbf{1}_{[-W,W]}(\xi)} \,d\xi
+
\int_{-W}^W \abs{\frac{1}{K}\rho _{K}(N,W;\xi)-\frac{N}{K}\mathbf{1}_{[-W,W]}(\xi)} \,d\xi
\\
&\qquad=
2W\left(\frac{N}{K}-\frac{1}{2W}\right)
+
\frac{2NW}{K}-\frac{1}{K}\sum_{k=0}^{K-1} \int_{-W}^W \abs{U_{k}(N,W;\xi)}^2\,d\xi
\\
&\qquad\leq \frac{2}{K} + 1-\frac{1}{K}\sum_{k=0}^{K-1} \lambda_k
\lesssim \frac{\log N}{K},
\end{align*}
thanks to Lemma \ref{lemma_sum_eig}. Now we estimate the broad brand leakage:
\begin{align*}
&\int_{I\setminus [-W,W]} \abs{\frac{1}{K}\rho _{K}(N,W;\xi)-\frac{1}{2W}\mathbf{1}_{[-W,W]}(\xi)} \,d\xi
=\int_{I\setminus [-W,W]} \frac{1}{K}\rho _{K}(N,W;\xi) \,d\xi
\\
&\qquad=\frac{1}{K}\sum_{k=0}^{K-1} (1-\lambda_k)
= 1 - \frac{1}{K}\sum_{k=0}^{K-1} \lambda_k,
\end{align*}
so the conclusion follows invoking again Lemma \ref{lemma_sum_eig}.
\section{Conclusion}

\label{sec_conclusion} In \cite[Section IV]{Thomson}, Thomson estimated $Bias(\ST)$ by using the approximation 
$\frac{1}{K}\rho _{K}(N,W,\cdot )\approx \frac{1}{2W}\mathbf{1}_{[-W,W]}$.
Besides supporting that reasoning, Theorem \ref{th_main} allows one to quantify the bias. Indeed, 
\begin{equation*}
\left\vert Bias(\ST(\xi ))\right\vert \leq \left\vert \left| S*\frac{1}{K}\rho_{K}(N,W,\cdot) - S*\frac{1}{2W}\mathbf{1}_{[-W,W]} \right|
\right\vert_\infty +\left\vert \left| S-S*\frac{1}{2W}\mathbf{1}_{[-W,W]} \right|
\right\vert_\infty
\end{equation*}
and, if $S$ is a bounded function, then Theorem \ref{th_main} implies
that 
\begin{equation*}
\left\vert \left| S*\frac{1}{K}\rho_{K}(N,W,\cdot)- S*\frac{1}{2W}\mathbf{1}_{[-W,W]} \right| \right\vert_\infty \lesssim \max_{\xi ^{\prime }\in \mathbb{R}}\left\vert S(\xi ^{\prime })\right\vert \frac{\log N}{K}.
\end{equation*}
The remaining term $\left\vert \left| S-S*\frac{1}{2W}\mathbf{1}_{[-W,W]}
\right| \right\vert_\infty $ can be bounded by assuming that $S$ is smooth. For
example, if, as in Thomson's work, $S$ is assumed to be analytic (and periodic), then $\left\vert \left| S-S*\frac{1}{2W}\mathbf{1}_{[-W,W]} \right| \right\vert
\lesssim W^{2}$, leading to the bias estimate:
\begin{equation}
\mathrm{Bias}(\ST(\xi ))\lesssim W^{2}+\frac{\log N}{K}.
\label{bias_estimate}
\end{equation}
On the other hand, for a slowly varying spectrum $S$, Thomson \cite{Thomson} argues that
\begin{equation}
\label{eq_var}
\mathrm{Var}\left( \ST(\xi )\right) \lesssim \frac{1}{K}\text{,}
\end{equation}
(see, \cite{WMP}, \cite{liro08} or \cite[Section 3.1.2]{HogLak} for precise expressions for the variance.) Given a 
number of available 
observations, the estimates in \eqref{bias_estimate} and \eqref{eq_var} show how much bias can be expected, in order to 
bring the variance down by a factor of $1/K$. This leads to a concrete estimate for the mean squared error
\begin{align}
\label{eq_MSE}
\mathrm{MSE}(\ST) = \mathbb{E}(S-\ST)^2
= \mathrm{Bias}(\ST)^2 + \mathrm{Var}(\ST)
\lesssim W^{4}+\frac{\log^2 N}{K^2} + \frac{1}{K},
\end{align}
that can be used to decide on the value of the bandwidth resolution parameter $W$.

We have thus obtained explicit bounds that allow us to quantify the
bias-variance trade-off in Thomson's multitaper method. Note that in the slowly varying regime, the error due to 
spectral leakage is largely dominated by the variance and therefore, in agreement with Thomson's analysis, the mean 
squared error is $\approx W^{4} + \frac{1}{K}$. In the case of more rapidly varying spectra, \eqref{eq_var} is no longer 
a valid approximation \cite{WMP, liro08} and the contribution of the spectral leakage to the mean squared error can be more 
significant.

\section{Appendix}
\subsection{Integral kernels}
The Toeplitz operator $H_{W}^{N}$ from \eqref{eq_toep_op} can be explicitly
described by the formula 
\begin{equation*}
H_{W}^{N}f(x)=\int_{I}f(y)K_{W}^{N}(x,y)dy\text{,}
\end{equation*}
where the kernel $K_{W}^{N}(x,y)$ is 
\begin{equation}
K_{W}^{N}(x,y)=\int_{[-W,W]}\mathbf{D}_{N}(x-z)\overline{\mathbf{D}_{N}(y-z)}dz\text{.}  \label{eq_kernel1}
\end{equation}

\subsection{An approximation lemma}
\begin{lemma}
\label{MainLemma} Let $f:I\rightarrow \mathbb{C}$ an integrable function, of
bounded variation, and supported on $I^{\circ }=(-1/2,1/2)$. For $N\geq 2$,
let 
\begin{equation*}
f\ast \left\vert \mathbf{D}_{N}\right\vert ^{2}(x)=\int_{I}f(y)\left\vert 
\mathbf{D}_{N}\left( x-y\right) \right\vert ^{2}dy,\qquad x\in I.
\end{equation*}
Then 
\begin{equation}
\left\Vert f-\frac{1}{N}f\ast \left\vert \mathbf{D}_{N}\right\vert
^{2}\right\Vert _{L^{1}(I)}\lesssim Var\left( f,I\right) \frac{\log N}{N}.
\end{equation}
\end{lemma}

\begin{remark}
In the above estimate, $Var(f,I)$ denotes the total variation of $f$ on $I$.
If $f=1_{[-W,W]}$, with $W\in (-1/2,1/2)$, then $Var\left( f,I\right) =2$
and the estimate reads 
\begin{equation*}
\left\Vert \mathbf{1}_{[-W,W]}-\frac{1}{N}\mathbf{1}_{[-W,W]}\ast \left\vert 
\mathbf{D}_{N}\right\vert ^{2}\right\Vert _{L^{1}(I)}\lesssim
\frac{\log N}{N}.
\end{equation*}
\end{remark}
\begin{proof}
By an approximation argument, we assume without loss of generality that $f$
is smooth (see for example \cite[Lemma 3.2]{AGR}). We also extend $f$ 
periodically to $\mathbb{R}$. Note that this extension is still smooth because $f|I$ is
supported on $I^{\circ }$.

\textbf{Step 1}. Since $f(x+h)-f(x)=\int_{0}^{1}f^{\prime }(th+x)h\,dt$, we
can use the periodicity of $f$ to estimate 
\begin{align*}
\left\Vert f(\cdot+h) - f \right\Vert _{L^{1}(I)} &\leq
\int_{0}^{1}\int_{-1/2}^{1/2}\left\vert f^{\prime }(t h +x)\right\vert
dx\left\vert h\right\vert dt =\int_{0}^{1}\int_{-1/2+th}^{1/2+th}\left\vert
f^{\prime }(x)\right\vert dx\left\vert h\right\vert dt \\
&\qquad =\int_{0}^{1}\int_{-1/2}^{1/2}\left\vert f^{\prime }(x)\right\vert
dx\left\vert h\right\vert dt=Var(f,I)\left\vert h\right\vert.
\end{align*}
Since $f$ is periodic, the previous estimate can be improved to: 
\begin{equation}
\left\Vert f(\cdot+h) - f\right\Vert _{L^{1}(I)}\lesssim Var(f,I)\left\vert
\sin(\pi h)\right\vert ,\qquad h\in \mathbb{R}.  \label{eq_f}
\end{equation}

\textbf{Step 2}. We use the notation $f^{N}:=f\ast \tfrac{1}{N}\left\vert 
\mathbf{D}_{N}\right\vert ^{2}$. By a change of variables and periodicity, 
\begin{equation*}
f(x)-f^{N}(x)=\frac{1}{N}\int_{-1/2}^{1/2}(f(x)-f(y+x))\left\vert
\mathbf{D}_{N}(-y)\right\vert ^{2}dy.
\end{equation*}
We can now finish the proof by resorting to \eqref{eq_f}: 
\begin{eqnarray*}
\left\Vert f-f^N\right\Vert _{L^{1}(I)} &\lesssim &Var(f,I)
\frac{1}{N}\int_{-1/2}^{1/2}\left\vert \sin (\pi y)\right\vert \left\vert
\mathbf{D}_{N}(y)\right\vert ^{2}dy \\
&\lesssim &Var(f,I)\frac{1}{N}\int_{0}^{1/2}\frac{\left\vert \sin (\pi
Ny)\right\vert }{\left\vert y\right\vert }dy \\
&\lesssim &Var(f,I)\frac{1}{N}\left[ 1+\int_{1}^{N/2}\frac{1}{\left\vert
y\right\vert }dy\right] \\
&\lesssim &Var(f,I)\frac{\log N}{N}\text{.}
\end{eqnarray*}
\end{proof}
\subsection{Proof of Lemma \ref{lemma_sum_eig}}
We first note from \eqref{eq_kernel1} that 
\begin{equation}
\trace\left( H_{W}^{N}\right)
=\int_{I}K_{W}^{N}(x,x)dx=\int_{[-W,W]}\int_{I}\left\vert
\mathbf{D}_{N}(x-y)\right\vert ^{2}dydx=2NW\text{,}  \label{eq_trace}
\end{equation}
since\ $\int_{I}$\ $\left\vert \mathbf{D}_{N}\right\vert ^{2}=N$. Moreover a
similar calculation gives 
\begin{equation*}
\trace\left( H_{W}^{N}\right) ^{2}=\int_{[-W,W]}\int_{I}\mathbf{1}_{_{[-W,W]}}(y)\left\vert \mathbf{D}_{N}(x-y)\right\vert ^{2}dydx\text{.}
\end{equation*}
Hence we can use Lemma \ref{MainLemma} to conclude that
\begin{align*}
\label{eq_trn}
&\trace\left[ \left( H_{W}^{N}\right) -\left( H_{W}^{N}\right) ^{2}\right]
=\int_{-W}^W\left[ N\mathbf{1}_{_{[-W,W]}}(x)-\left( \mathbf{1}_{_{[-W,W]}}\ast
\left\vert \mathbf{D}_{N}\right\vert ^{2}\right) (x)\right] dx
\\
&\qquad \leq
\int_{I}\abs{ N\mathbf{1}_{_{[-W,W]}}(x)-\left( \mathbf{1}_{_{[-W,W]}}\ast
\left\vert \mathbf{D}_{N}\right\vert ^{2}\right) (x)} dx
\leq C \log N,
\end{align*}
for some constant $C$. Using this bound, we estimate:
\begin{align*}
C\log N &\geq \sum_{k=0}^{N-1} \lambda_k(1-\lambda_k)
=\sum_{k=0}^{K-1} \lambda_k(1-\lambda_k)
+\sum_{k=K}^{N-1} \lambda_k(1-\lambda_k)
\\
&\geq \lambda_{K-1} \sum_{k=0}^{K-1} (1-\lambda_k)
+ (1-\lambda_{K-1})\sum_{k=K}^{N-1} \lambda_k
\\
&= \lambda_{K-1} K - \lambda_{K-1}\sum_{k=0}^{K-1} \lambda_k
+ (1-\lambda_{K-1})(2NW-\sum_{k=0}^{K-1} \lambda_k)
\\
&= \lambda_{K-1} K
+ 2NW(1-\lambda_{K-1}) - \sum_{k=0}^{K-1} \lambda_k
\\
&= 2NW - \sum_{k=0}^{K-1} \lambda_k
+
\lambda_{K-1} (K-2NW)
\\
&\geq K - \sum_{k=0}^{K-1} \lambda_k
-1.
\end{align*}
Hence, $K - \sum_{k=0}^{K-1} \lambda_k \leq C\log N + 1$. On the other
hand $\sum_{k=0}^{K-1} \lambda_k - K \leq 2NW - K \leq 1$. Therefore,
since $N \geq 2$, $\abs{K - \sum_{k=0}^{K-1} \lambda_k} \lesssim \log N$ and the conclusion follows.


\begin{thebibliography}{99}
\bibitem{AB} L. D. Abreu, A. S. Bandeira, \emph{Landau's necessary density
conditions for the Hankel transform},\ \ J. Funct. Anal. 162 (2012),
1845-1866.

\bibitem{AGR} L. D. Abreu, K. Gr\"{o}chenig, J. L. Romero, \emph{On
accumulated spectrograms},\emph{\ }Trans. Amer. Math. Soc., 368 (2016),
3629-3649.

\bibitem{BB} M. Bayram, R. G. Baraniuk \emph{Multiple window time-varying
spectrum estimation},\emph{\ }In Nonlinear and Nonstationary Signal
Processing (Cambridge, 1998), pages 292-316. Cambridge Univ. Press,
Cambridge, 2000.

\bibitem{Benedetto} J. J. Benedetto, \emph{Harmonic analysis and spectral
estimation}, J. Math. Anal. Appl., \textbf{91} (1983), 444-509.

\bibitem{Science1997} G. Bond, W. Showers, M. Cheseby, R. Lotti, P. Almasi,
P. deMenocal, P. Priore, H. Cullen, I. Hajdas, G. Bonani, \emph{A Pervasive
Millennial-Scale Cycle in North Atlantic Holocene and Glacial Climates},
Science (1997), \textbf{278}, 1257-1266.

\bibitem{SimmonsGeoph2008} F. A. Dahlen, F. J. Simons, \emph{Spectral
estimation on a sphere in geophysics and cosmology}, Geophys. J. Int. (2008) 
\textbf{174}, 774-807.

\bibitem{Neurosciences} A. Delorme, S. Makeig, \emph{EEGLAB: an open source
toolbox for analysis of single-trial EEG dynamics including independent
component analysis} - J. Neurosci. Methods, (2004) \textbf{134}, 9-21.

\bibitem{Cognitive} S. Haykin, D. J. Thomson, J. H. Reed, \emph{Spectrum
sensing for cognitive radio}, Proc. IEEE, (2009), \textbf{97}, 849 - 877.

\bibitem{PNAS2012} C. Harig, F. J. Simons, \emph{Mapping Greenland's mass
loss in space and time}, Proc. Natl. Acad. Sci. USA, (2012), \textbf{109},19934--19937.

\bibitem{Haykin} S. Haykin, \emph{Cognitive radio: brain-empowered wireless
communications, } IEEE Journal on Selected Areas in Communications,
\textbf{23} (2), 201-220, 2005.

\bibitem{HogLak} J. A. Hogan, J. D. Lakey, \emph{Duration and Bandwith
Limiting. Prolate Functions, Sampling, and Applications}, Applied and
Numerical Harmonic Analysis, Birkh\"{a}user/Springer, New York, 2012,
xvii+258pp.

\bibitem{la67-1}
H.~J. {L}andau.
\newblock {S}ampling, data transmission, and the {N}yquist rate.
\newblock {\em Proc. IEEE}, 55(10):1701--1706, {O}ctober 1967.

\bibitem{la75-1}
H.~J. {L}andau.
\newblock {O}n {S}zeg{\"o}'s eigenvalue distribution theorem and
  non-{H}ermitian kernels.
\newblock {\em J. Anal. Math.}, 28:335--357, 1975.

\bibitem{lapo61}
H.~J. {L}andau and H.~O. {P}ollak.
\newblock {P}rolate spheroidal wave functions, {F}ourier analysis and
uncertainty {I}{I}.
\newblock {\em Bell System Tech. J.}, 40:65--84, 1961.

\bibitem{liro08}
K. S. Lii and M. Rosenblatt, 
\newblock {P}rolate spheroidal spectral estimates.
\newblock {\em Stat. Probab. Lett.}, 78 (11), 1339-1348, 2008.

\bibitem{OW} S. C. Olhede, A. T. Walden,
\emph{Generalized Morse wavelets},
IEEE Trans. Signal Process, 50 (11), 2661-2670, 2002.

\bibitem{PW} D. B. Percival, A. T. Walden, \emph{Spectral Analysis for
Physical Applications, Multitaper and Conventional Univariate Techniques.
Cambridge,} 1993, Cambridge University Press.

\bibitem{ACHA2014} A. Plattner, F. J. Simons, \emph{Spatiospectral
concentration of vector fields on a sphere}, Appl. Comp. Harm. Anal. 36 (1),
(2014) 1-22.

\bibitem{SIAMRev} F.~J. Simons, F.~A. Dahlen, and M.~A. Wieczorek.
\newblock Spatiospectral concentration on a sphere.
\newblock {\em SIAM Rev.}, 48(3):504--536 (electronic), 2006.

\bibitem{Slepian} D. Slepian, \emph{Some comments on Fourier analysis,
uncertainty and modeling},\emph{\ }SIAM Rev. \textbf{25} (1983) 379-393.

\bibitem{SlepianIV} D. Slepian, \emph{Prolate spheroidal wave
functions, Fourier analysis and uncertainty-IV: Extensions to many
dimensions; generalized prolate spheroidal functions \ }Bell Syst. Tech. J.,
(1964), 3009-3057.

\bibitem{posl61}
D.~{S}lepian and H.~O. {P}ollak.
\newblock {P}rolate {S}pheroidal {W}ave {F}unctions, {F}ourier {A}nalysis and
{U}ncertainty {I}.
\newblock {\em {I}. {B}ell {S}yst. {T}ech.{J}.}, 40(1):43--63, 1961. 


\bibitem{Thomson} D. J. Thomson, \emph{Spectrum estimation and harmonic
analysis},\emph{\ }Proc.IEEE, \textbf{70}, (1982) 1055-1095.\textsl{\ }

\bibitem{Thomson2} D. J. Thomson, \emph{Multitaper Analysis of Nonstationary
and Nonlinear Time Series Data}, Nonlinear and Nonstationary Signal
Processing, Cambridge University Press, (2000).

\bibitem{WMP} A. T. Walden, E. J. McCoy, D. B. Percival, \emph{The variance
of multitaper spectrum estimates for real Gaussian processes}. IEEE Trans.
Signal Process, \textbf{42} (1994), 479-482.

\bibitem{zhsi13} Z. Zhao, A. Singer, \emph{Fourier-Bessel Rotational
Invariant Eigenimages}, The Journal of the Optical Society of America A, 30
(5), pp. 871--877 (2013).
\end{thebibliography}
\end{document}